\documentclass{amsart}
\usepackage{amssymb}

\newtheorem{theorem}{Theorem}

\theoremstyle{definition}

\theoremstyle{remark}
\newtheorem*{remark}{Remark}
\newtheorem*{note}{Note}

\numberwithin{equation}{section}

\begin{document}
\rightline{\scriptsize CRM-2442}

\title[Multiple Bailey, Rogers and Dougall sums]
{On certain multiple Bailey, Rogers and \\ Dougall type summation
formulas}

\author{J. F. van Diejen}
\address{Centre de Recherches Math\'ematiques,
Universit\'e de Montr\'eal,  C.P. 6128, succursale Centre-ville,
Montr\'eal (Qu\'ebec), H3C 3J7 Canada}

\thanks{Work supported in part by
the Natural Sciences and Engineering Research Council (NSERC) of Canada.}


\date{December 1996}


\begin{abstract}
A multidimensional generalization of Bailey's very-well-poised
bilateral basic hypergeometric ${}_6\psi_6$ summation formula and its
Dougall type ${}_5H_5$ hypergeometric degeneration for $q\rightarrow 1$
is studied.
The multiple Bailey sum amounts to an extension corresponding
to the case of a nonreduced root system of certain summation identities
associated to the reduced root systems that were recently conjectured by
Aomoto and Ito and proved by Macdonald.
By truncation, we obtain multidimensional analogues of
the very-well-poised unilateral (basic) hypergeometric Rogers
${}_6\phi_5$ and Dougall
${}_5F_4$ sums (both nonterminating and terminating).
The terminating sums may be used to arrive at product formulas
for the norms of recently introduced ($q$-)Racah polynomials
in several variables.
\end{abstract}

\maketitle

\section{Introduction}
In this paper certain multidimensional generalizations are studied
of the summation formulas
\begin{subequations}
\begin{eqnarray}\label{bai-sum}
\lefteqn{\sum_{\lambda \in \mathbb{Z}}\,
\frac{\prod_{r=1}^{4} ( q^{1+g_r+z+\lambda}, q^{1+g_r-z-\lambda};q)_\infty}
     { (q^{1+2z+2\lambda}, q^{1-2z-2\lambda};q)_\infty } } && \\
&& \makebox[10em]{} =
     \frac{(q;q)_\infty \prod_{1\leq r < s\leq 4} (q^{1+g_r+g_s};q)_\infty}
        {(q^{1+g_1+g_2+g_3+g_4};q)_\infty}  \nonumber
\end{eqnarray}
(with $2z\not\in (\mathbb{Z}+\frac{2\pi }{i \log q}\mathbb{Z})$) and
\begin{eqnarray}\label{doug-sum}
\lefteqn{\sum_{\lambda \in \mathbb{Z}}\,
\frac{\Gamma (1+2z+2\lambda ) \Gamma (1-2z-2\lambda )}
     {\prod_{r=1}^{4} \Gamma (1+g_r+z+\lambda )
                      \Gamma (1+g_r-z-\lambda )}} && \\
& & \makebox[10em]{} =\frac{\Gamma (1+g_1+g_2+g_3+g_4)}
     {\prod_{1\leq r < s\leq 4} \Gamma (1+g_r+g_s)} \nonumber
\end{eqnarray}
\end{subequations}
(with $2z\not\in \mathbb{Z}$),
where it is assumed that
$0<q<1$ and $\text{Re}(1+g_1+g_2+g_3+g_4)>0$.
(For conventions regarding the notation,
we refer to the remark at the end of this introduction.)
The conditions on $z$, $q$ and on the parameters $g_1,g_2,g_3,g_4$
guarantee that the terms on the l.h.s. are finite and, moreover, that the
series converge in absolute value.
The sum in \eqref{bai-sum}
is (for generic $z$) equivalent to Bailey's very-well-poised
bilateral ${}_6\psi_6$ sum \cite{bai:series,gas-rah:basic}
\begin{subequations}
\begin{eqnarray}\label{bai-ser}
\lefteqn{
\makebox[2em]{}\sum_{\lambda \in \mathbb{Z}}
q^{(1+g_1+g_2+g_3+g_4)\lambda}
\Big( \frac{1-q^{2z+2\lambda}}
           {1-q^{2z}} \Bigr)
\prod_{1\leq r\leq 4}
\frac{(q^{z-g_r};q)_\lambda}
     {(q^{1+g_r+z};q)_\lambda} } &&  \\
&&={}_6\psi_6 \Bigl(
\begin{array}{c}
q^{1+z},-q^{1+z}, q^{z-g_1},q^{z-g_2},q^{z-g_3},q^{z-g_4} \\
q^z, -q^z, q^{1+g_1+z},q^{1+g_2+z},q^{1+g_3+z},q^{1+g_4+z}
\end{array} ;q, q^{1+g_1+g_2+g_3+g_4} \Bigr)  \nonumber \\
&& =\frac{ (q^{1+2z}, q^{1-2z};q)_\infty }
        {\prod_{r=1}^{4} (q^{1+g_r+z}, q^{1+g_r-z};q)_\infty}
\frac{ (q ;q)_\infty \prod_{1\leq r < s\leq 4} (q^{1+g_r+g_s};q)_\infty}
     {(q^{1+g_1+g_2+g_3+g_4} ;q)_\infty}  \nonumber
\end{eqnarray}
($0<q<1$) and the sum in \eqref{doug-sum} corresponds to the degenerate case
$q\rightarrow 1$, which amounts to Dougall's
very-well-poised bilateral ${}_5H_5$ sum
\cite{dou:vandermonde,bai:series}
\begin{eqnarray}\label{doug-ser}
\lefteqn{
\sum_{\lambda \in \mathbb{Z}}\:
 \Bigl( 1+\frac{\lambda}{z}\Bigr)
\prod_{1\leq r\leq 4}\frac{(z-g_r)_\lambda}
                          {(1+z+g_r)_\lambda} } && \\
&& =
{}_5H_5 \left(
\begin{array}{c}
1+z, z-g_1,z-g_2,z-g_3,z-g_4 \\
z, 1+z+g_1,1+z+g_2,1+z+g_3+1+z+g_4
\end{array}; 1 \right) \nonumber \\
&& =
\frac{\prod_{r=1}^{4} \Gamma (1+g_r+z) \Gamma (1+g_r-z)}
     {\Gamma (1+2z) \Gamma (1-2z)}
\frac{\Gamma(1+g_1+g_2+g_3+g_4)}
     {\prod_{1\leq r < s\leq 4} \Gamma (1+g_r+g_s)} \nonumber
\end{eqnarray}
\end{subequations}
(with the same convergence condition on the parameters
$\text{Re}\: (1+g_1+g_2+g_3+g_4)>0$).
In order for the terms of the series in \eqref{bai-ser}, \eqref{doug-ser}
to be finite, the variable $z$ should be chosen such that it
is nonzero (modulo $\frac{2\pi }{i\log q}$) and such that
$g_r\pm z$ is not a negative integer (modulo $\frac{2\pi }{i\log q}$).
The identities in \eqref{bai-sum} and \eqref{doug-sum}
pass over into the Bailey and Dougall sums in
\eqref{bai-ser}, \eqref{doug-ser}
after division by the middle term corresponding to
$\lambda =0$ (which is nonzero with these restrictions on $z$)
and rewriting of the resulting l.h.s. with
the aid of the relations
$(a;q)_\infty /(aq^\lambda ;q)_\infty = (a;q)_\lambda$ and
$\Gamma (a+\lambda)/\Gamma (a) =(a)_\lambda$ as well as
the reflection properties
$(a;q)_{\lambda}(a^{-1}q;q)_{-\lambda}=
(-a)^{\lambda}q^{\lambda (\lambda -1)/2}$
and $(a)_{\lambda}(1-a)_{-\lambda} = (-1)^\lambda$, respectively.

The plan of the paper reads as follows.
In Section~\ref{sec2} a multidimensional version of
the sums in \eqref{bai-sum} and \eqref{doug-sum} is discussed
(Theorem~\ref{sumMthm}).
For $0<q<1$ the sum under consideration amounts
to a generalization corresponding
to the case of a nonreduced root system
of summation formulas appearing
in recent work of Macdonald \cite{mac:formal}
associated to the (affine) reduced root systems.
Alternative representations lead to certain multidimensional
analogues of the Bailey ${}_6\psi_6$ sum \eqref{bai-ser} and the
Dougall ${}_5H_5$ sum \eqref{doug-ser},
as well as to multiple summation formulas
of the type studied by Aomoto and Ito \cite{aom:product,ito:theta}.
Section~\ref{sec3} describes the specialization to nonterminating and
terminating
multidimensional versions of Rogers' very-well-poised ${}_6\phi_5$ sum
\cite{rog:third,gas-rah:basic}
and Dougall's very-well-poised ${}_5F_4$ sum
\cite{dou:vandermonde,gas-rah:basic} (Theorem~\ref{sumRDthm} and
\ref{finsumRDthm}).
The resulting terminating sums (of Theorem~\ref{finsumRDthm}) may be
used to arrive at product formulas for the norms of recently introduced
($q$-)Racah polynomials in several variables \cite{die-sto:multivariable}
that generalize the well-known one-variable ($q$-)Racah polynomials of
Askey and Wilson \cite{ask-wil:set,wil:some,gas-rah:basic}.
The technicalities pertaining to the proof for the multiple bilateral
summation formulas of Section~\ref{sec2} are relegated to Section~\ref{sec4}
and an appendix at the end of the paper (in which
the convergence of the series is demonstrated).
The proof in Section~\ref{sec4} is based on a recurrence relation for the
generalized
Macdonald type sum (of Theorem~\ref{sumMthm}), which is derived
using a technique very similar to that employed by
Gustafson in his proof of the Selberg type multivariable
Askey-Wilson integral studied in \cite{gus:generalization}.

The sums considered in this paper are not the
only/first possible (nontrivial)
multidimensional generalizations of the ${}_5F_4$, ${}_6\phi_5$,
${}_5H_5$ and ${}_6\psi_6$ summation formulas.
An important class of very-well-poised (basic) hypergeometric
summation formulas associated with the (special)
unitary group $(S)U(n)$ (type $A$ root system) can e.g. be
found in the works of Holman \cite{hol:summation} (${}_5F_4$ type),
Milne \cite{mil:q-analog,mil:basic} (${}_6\phi_5$ type) and
Gustafson \cite{gus:multilateral} (${}_5H_5$ and ${}_6\psi_6$ type).
Gustafson moreover generalized his
$(S)U(n)$ type multiple ${}_5H_5$ and ${}_6\psi_6$
sums of \cite{gus:multilateral} to the case of an arbitrary
classical simple Lie group \cite{gus:macdonald}.
For the symplectic group $Sp(n)$ (type $C$ root system)
truncated versions of Gustafson's ${}_6\psi_6$ sum
giving rise to terminating and nonterminating
multiple ${}_6\phi_5$ summation formulas
were discussed by Lilly and Milne \cite{lil-mil:bailey}.
Very recently, still other multidimensional versions of the
${}_6\psi_6$ and ${}_6\phi_5$ summation formulas associated
with the type $C$ root system were
presented by Schlosser \cite{sch:summation}.

All these multidimensional generalizations of
the ${}_5F_4$, ${}_6\phi_5$,
${}_5H_5$ and ${}_6\psi_6$ summation formulas
occurring in the literature
are of a different type than those studied below.
We will, however, have the opportunity in Section~\ref{sec4} to
employ Gustafson's multiple
${}_5H_5$ and ${}_6\psi_6$ summation
formulas from \cite{gus:macdonald}
for the symplectic group $Sp(n)$ (type $C$ root system),
when deriving the recurrence relation that leads to the proof of
the Macdonald type sum given by Theorem~\ref{sumMthm}.

\begin{remark}
We have adopted the following (standard)
conventions regarding
the notation of $q$-shifted factorials and Pochhammer symbols,
respectively \cite{gas-rah:basic}
\begin{subequations}
\begin{equation}
(a;q)_m  \equiv \left\{
\begin{array}{lll}
1 & \text{for} & m=0 \\
(1-a)(1-aq)\cdots (1-aq^{m-1}) & \text{for} & m=1,2,3,\ldots \\
\frac{1}{(1-aq^{-1})(1-aq^{-2}) \cdots (1-aq^m)} & \text{for} &
m= -1,-2,-3,\ldots
\end{array}\right.
\end{equation}
and
\begin{equation}
(a)_m  \equiv \left\{
\begin{array}{lll}
1 & \text{for} & m=0 \\
a(a+1)\cdots (a+m-1) & \text{for} & m=1,2,3,\ldots \\
\frac{1}{(a-1)(a-2) \cdots (a+m)} & \text{for} & m= -1,-2,-3,\ldots
\end{array}\right.
\end{equation}
\end{subequations}
(where for negative $m$ it is assumed that the value of $a$ is such that
the denominator does not vanish).
One has
\begin{equation}
(a;q)_m = \frac{(a;q)_\infty}{(aq^m;q)_\infty}
\;\;\;\;\;\;\;\;\; \text{and} \;\;\;\;\;\;\;\;\;
(a)_m = \frac{\Gamma (a+m)}{\Gamma (a)},
\end{equation}
where $\Gamma (\cdot)$ represents the gamma function and
\begin{equation}
(a;q)_\infty \equiv \prod_{k=0}^\infty (1-aq^k)
\end{equation}
(here it is assumed that $|q|<1$).
We furthermore use the abbreviation
\begin{subequations}
\begin{eqnarray}
(a_1,\ldots ,a_p;q)_m &\equiv & (a_1;q)_m\cdots (a_p;q)_m ,\\
(a_1,\ldots ,a_p)_m &\equiv & (a_1)_m\cdots (a_p)_m.
\end{eqnarray}
\end{subequations}
One-variable bilateral
(basic) hypergeometric series and their unilateral counterparts
are denoted by
\begin{subequations}
\begin{eqnarray}
{}_p\psi_p \left(
\begin{array}{c}
a_1,\ldots ,a_p \\
b_1,\ldots ,b_p
\end{array} ; q, \zeta \right)
&\equiv&
\sum_{\lambda\in \mathbb{Z}}
\frac{(a_1,\ldots ,a_p;q)_\lambda}
     {(b_1,\ldots ,b_p;q)_\lambda}\: \zeta^\lambda , \\
{}_pH_p \left(
\begin{array}{c}
a_1,\ldots ,a_p \\
b_1,\ldots ,b_p
\end{array} ; \zeta \right)
&\equiv&
\sum_{m\in \mathbb{Z}}
\frac{(a_1,\ldots ,a_p)_\lambda}
     {(b_1,\ldots ,b_p)_\lambda}\: \zeta^\lambda
\end{eqnarray}
\end{subequations}
and
\begin{subequations}
\begin{eqnarray}
{}_{p}\phi_{p-1} \left(
\begin{array}{c}
a_1,\ldots ,a_{p} \\
b_1,\ldots ,b_{p-1}
\end{array} ; q, \zeta \right)
&\equiv&
\sum_{\lambda\in \mathbb{N}}
\frac{(a_1,\ldots ,a_{p};q)_\lambda}
     {(b_1,\ldots ,b_{p-1},q;q)_\lambda}\: \zeta^\lambda , \\
{}_{p}F_{p-1} \left(
\begin{array}{c}
a_1,\ldots ,a_{p} \\
b_1,\ldots ,b_{p-1}
\end{array} ; \zeta \right)
&\equiv &
\sum_{m\in \mathbb{N}}
\frac{(a_1,\ldots ,a_{p})_\lambda}
     {(b_1,\ldots ,b_{p-1},1,q)_\lambda}\: \zeta^\lambda
\end{eqnarray}
\end{subequations}
(where $\mathbb{N}$ includes the number zero).
\end{remark}

\section{Multiple analogues of very-well-poised bilateral\\
(basic) hypergeometric series}\label{sec2}
\begin{note}
In this section it is always assumed that the nome
$q$ lies in the open interval $]0,1[$ (corresponding to
the case of basic hypergeometric series)
or that it degenerates to $q=1$ (corresponding to
the case of ordinary hypergeometric series).
\end{note}

\subsection{Notation}
We will first set up some notational preliminaries.
The reader may wish to skip this part at first reading and
refer back to it when needed.

To describe the multiple summation formulas below
it is convenient to introduce the functions
$\mathcal{C}_{+,q}(\mathbf{x})$,
$\mathcal{C}_{-,q}(\mathbf{x})$,
$\hat{\mathcal{C}}_{+,q}(\mathbf{x})$ and
$\hat{\mathcal{C}}_{+,q}(\mathbf{x})$ given
for $0< q <1$ by
\begin{subequations}
\begin{eqnarray}\label{C+q}
\mathcal{C}_{+,q} (\mathbf{x})&=&
\prod_{1\leq j<k\leq n}
\frac{(q^{1+x_j+x_k}, q^{1+x_j-x_k};q)_\infty}
{(q^{1+g+x_j+x_k},q^{1+g+x_j-x_k};q)_\infty} \\
&&\times \prod_{1\leq j\leq n}
\frac{(q^{1+2x_j};q)_\infty }
{\prod_{r=1}^{4} ( q^{1+g_r+x_j};q)_\infty} , \nonumber \\
\mathcal{C}_{-,q} (\mathbf{x})&=&
q^{-\sum_{j=1}^n (1+2\hat{\rho}_j) x_j}
\prod_{1\leq j<k\leq n}
\frac{(q^{-g+x_j+x_k},q^{-g+x_j-x_k};q)_\infty}
{(q^{x_j+x_k}, q^{x_j-x_k};q)_\infty} \\
&&\times \prod_{1\leq j\leq n}
\frac{\prod_{r=1}^{4} ( q^{-g_r+x_j};q)_\infty}
{(q^{2x_j};q)_\infty }  ,\nonumber \\
\hat{\mathcal{C}}_{+,q} (\mathbf{x})&=&
\prod_{1\leq j<k\leq n}
\frac{(q^{1+x_j+x_k}, q^{1+x_j-x_k};q)_\infty}
{(q^{1+g+x_j+x_k},q^{1+g+x_j-x_k};q)_\infty} \\
&&\times \prod_{1\leq j\leq n}
\frac{(q^{1+2x_j};q)_\infty }
{\prod_{r=1}^{4} ( q^{1+\hat{g}_r+x_j};q)_\infty}  ,\nonumber
\end{eqnarray}
\begin{eqnarray}
\hat{\mathcal{C}}_{-,q} (\mathbf{x})&=&
\prod_{1\leq j<k\leq n}
\frac{(q^{1-g+x_j+x_k},q^{1-g+x_j-x_k};q)_\infty}
{(q^{1+x_j+x_k}, q^{1+x_j-x_k};q)_\infty} \\
&&\times \prod_{1\leq j\leq n}
\frac{\prod_{r=1}^{4} ( q^{1-\hat{g}_r+x_j};q)_\infty}
{(q^{1+2x_j};q)_\infty } , \nonumber
\end{eqnarray}
\end{subequations}
and for the degenerate case $q=1$ by
\begin{subequations}
\begin{eqnarray}\label{C+1}
\mathcal{C}_{+,1} (\mathbf{x})&=&
\prod_{1\leq j<k\leq n}
\frac{\Gamma (1+g+x_j+x_k)\, \Gamma(1+g+x_j-x_k)}
{\Gamma (1+x_j+x_k)\,\Gamma (1+x_j-x_k)} \\
&&\times \prod_{1\leq j\leq n}
\frac{\prod_{r=1}^{4} \Gamma (1+g_r+x_j)}
{\Gamma(1+2x_j)}  ,\nonumber \\
\mathcal{C}_{-,1} (\mathbf{x})&=&
\prod_{1\leq j<k\leq n}
\frac{\Gamma (x_j+x_k)\,\Gamma(x_j-x_k)}
{\Gamma(-g+x_j+x_k)\,\Gamma(-g+x_j-x_k)} \\
&&\times \prod_{1\leq j\leq n}
\frac{\Gamma (2x_j) }
{\prod_{r=1}^{4}  \Gamma (-g_r+x_j)}   ,\nonumber \\
\hat{\mathcal{C}}_{+,1} (\mathbf{x})&=&
\prod_{1\leq j<k\leq n}
\frac{\Gamma (1+g+x_j+x_k)\,\Gamma(1+g+x_j-x_k)}
{\Gamma(1+x_j+x_k)\, \Gamma(1+x_j-x_k)}\\
&&\times \prod_{1\leq j\leq n}
\frac{\prod_{r=1}^{4} \Gamma ( 1+\hat{g}_r+x_j)}
{\Gamma(1+2x_j) } ,\nonumber \\
\hat{\mathcal{C}}_{-,1} (\mathbf{x})&=&
\prod_{1\leq j<k\leq n}
\frac{\Gamma (1+x_j+x_k)\,\Gamma(1+x_j-x_k)}
{\Gamma(1-g+x_j+x_k)\, \Gamma(1-g+x_j-x_k)} \\
&&\times \prod_{1\leq j\leq n}
\frac{\Gamma(1+2x_j) }
{\prod_{r=1}^{4} \Gamma ( 1-\hat{g}_r+x_j)} .\nonumber
\end{eqnarray}
\end{subequations}
The parameters $\hat{g}_r$, $r=1,2,3,4$, are related to
the parameters $g_r$, $r=1,2,3,4$, by the transformation
\begin{equation}
\left(
\begin{array}{r}
\hat{g}_a \\ \hat{g}_b \\ \hat{g}_c \\ \hat{g}_d
\end{array} \right) =
\frac{1}{2}
\left(
\begin{array}{rrrr}
1 & 1 & 1 & 1 \\
1 & 1 & -1 & -1 \\
1 & -1 & 1 & -1 \\
1 & -1 &-1 & 1
\end{array} \right)
\left(
\begin{array}{r}
g_a \\ g_b \\ g_c \\g_d
\end{array} \right)  ,
\end{equation}
where $a$, $b$, $c$ and $d$ denote
a (fixed but otherwise arbitrary) permutation of the indices
$1$, $2$, $3$ and $4$ (i.e., one has that
$\{ a,b,c,d \} = \{ 1,2,3,4 \}$). We may in fact set $a=1$, $b=2$,
$c=3$ and $d=4$ without loss of generality, but here we have preferred
not to fix such choice explicitly in order to reflect
in our notation the invariance of the construction with respect to
permutations of the parameters $g_1,g_2,g_3 ,g_4$.

We will furthermore employ the vectors $\rho$ and $\hat{\rho}$
with components given by
\begin{equation}\label{rho}
\rho_j = (n-j)g +g_a\;\;\;\;\;\; \text{and}\;\;\;\;\;\;
\hat{\rho}_j=(n-j)g+\hat{g}_a
\end{equation}
($j=1,\ldots ,n$) and the Jacobi theta function
\begin{equation}\label{theta}
\theta (\zeta ) = \sum_{m=-\infty}^\infty (-1)^m q^{m(m-1)/2} \zeta^m .
\end{equation}
This theta function satisfies the quasi-periodicity relation
\begin{equation}\label{quasi}
\theta (q\zeta ) = -\zeta^{-1} \theta (\zeta)
\end{equation}
and admits the product representation
\begin{equation}\label{prod}
\theta (\zeta ) = (q,\zeta , q\zeta^{-1};q)_\infty .
\end{equation}
The equality of the r.h.s.
of \eqref{theta} and \eqref{prod} hinges on a
classic bilateral
summation formula known as the Jacobi triple product identity
(see e.g. \cite{gas-rah:basic}). The $q=1$ counterpart of
$\theta (q^z )$ is given by the sine function $\sin (\pi z)$
and the analogue of the product formula boils in this
degenerate situation down to the
reflection relation for the gamma function
\begin{equation}\label{refg}
\sin (\pi z ) = \frac{\pi}{\Gamma (z)\, \Gamma (1-z)}.
\end{equation}
The corresponding period lattice $\Omega_q\subset\mathbb{C}$
is given by
\begin{equation}\label{lattice}
\Omega_q = \left\{
\begin{array}{ll}
\mathbb{Z} +\frac{2\pi }{i\log (q)}\mathbb{Z} & \text{for}\;\; 0<q<1 \\ [1ex]
\mathbb{Z} & \text{for} \;\; q=1 .
\end{array} \right.
\end{equation}

\subsection{A Macdonald type sum}\label{sec2.2}
The following theorem describes a multidimensional generalization
of the summation formulas \eqref{bai-sum} (when $0<q<1$) and
\eqref{doug-sum} (when $q=1$) and reduces to these
formulas for $n=1$.
\begin{theorem}\label{sumMthm}
Let $0< q\leq 1$. For parameters subject to the condition
\begin{equation}\label{parco}
\text{Re}\, (1+2(n-j)g+g_1+g_2+g_3+g_4) > 0
\end{equation}
(with $j=1,\ldots ,n$), one has that
\begin{equation}\label{sumM}
\sum_{\lambda \in \mathbb{Z}^n}
\frac{1}
     {\mathcal{C}_{+,q}(\mathbf{z}+\lambda )\,
      \mathcal{C}_{+,q}(-\mathbf{z}-\lambda ) }
=
\frac{\hat{\mathcal{C}}_{-,q}(\hat{\rho})}
     {\hat{\mathcal{C}}_{+,q}(\hat{\rho})}
\end{equation}
and the series on the l.h.s. converges in absolute value.
Here it is assumed
that $\mathbf{z}\in\mathbb{C}^n$ with
the combinations
$z_j+z_k$, $z_j-z_k$ ($1\leq j<k\leq n$) and $2z_j$ ($1\leq j\leq n$)
being nonzero modulo
the lattice $\Omega_q$ \eqref{lattice}.
\end{theorem}
The genericity
condition on the components of the vector $\mathbf{z}$ ensures that
the denominators of the terms on the l.h.s. do not vanish;
the proof that the series under consideration converges in
absolute value can be found in the appendix at the end of the paper.
The evaluation of the sum hinges on a recurrence relation
that is derived by means
of a technique due to Gustafson \cite{gus:generalization},
who used it
to evaluate a Selberg type multivariable generalization
of the Askey-Wilson integral (see \cite{gas-rah:basic} and
references therein for
a discussion of the Askey-Wilson integral). The details of this
derivation, leading to a proof for the fact
that the value of the sum is given by the $\mathbf{z}$-independent
constant on the r.h.s. of \eqref{sumM}, are relegated to Section~\ref{sec4}.

A basic hypergeometric summation formula
closely related to that described by the theorem
has been derived recently by Macdonald using the
properties of affine Weyl groups and root systems
\cite{mac:formal}.
More specifically, the sum considered by
Macdonald is associated to an arbitrary reduced root
system and Theorem~\ref{sumMthm} (with $0<q<1$) may be viewed as the
extension to the case of a nonreduced root system.
For special values of the parameters $g_1,g_2,g_3 ,g_4$,
the basic hypergeometric sum in \eqref{sumM} reduces
to the Macdonald sums related to the reduced root systems of classical type
corresponding to the $B$, $C$ and $D$ series.

The evaluation constant
on the r.h.s. of \eqref{sumM} given by
\begin{subequations}
\begin{eqnarray}\label{normMq}
\makebox[1em]{} \frac{\hat{\mathcal{C}}_{-,q}(\hat{\rho})}
                     {\hat{\mathcal{C}}_{+,q}(\hat{\rho})}
\!\!\!\!\! &=&\!\!\!\!\!
\prod_{1\leq j<k\leq n}
\frac{(q^{1+g+\hat{\rho}_j+\hat{\rho}_k},
       q^{1+g+\hat{\rho}_j-\hat{\rho}_k};q)_\infty}
     {(q^{1+\hat{\rho}_j+\hat{\rho}_k},
       q^{1+\hat{\rho}_j-\hat{\rho}_k};q)_\infty} \\
&&\makebox[4em]{} \times
\frac{(q^{1-g+\hat{\rho}_j+\hat{\rho}_k},
       q^{1-g+\hat{\rho}_j-\hat{\rho}_k};q)_\infty}
     {(q^{1+\rho_j+\rho_k},
       q^{1+\hat{\rho}_j-\hat{\rho}_k};q)_\infty} \nonumber \\
&&\times \prod_{1\leq j\leq n}
\frac{\prod_{r=1}^{4}
      (q^{1+\hat{g}_r+\hat{\rho}_j}, q^{1-\hat{g}_r+\hat{\rho}_j};q)_\infty }
     {(q^{1+2\hat{\rho}_j};q)_\infty^2}  \nonumber
\end{eqnarray}
for $0<q<1$ and by
\begin{eqnarray}\label{normM1}
\makebox[1em]{} \frac{\hat{\mathcal{C}}_{-,1}(\hat{\rho})}
                     {\hat{\mathcal{C}}_{+,1}(\hat{\rho})}
\!\!\!\!\!&=&\!\!\!\!\!
\prod_{1\leq j<k\leq n}
\frac{\Gamma(1+\hat{\rho}_j+\hat{\rho}_k)\Gamma(1+\hat{\rho}_j-\hat{\rho}_k)}
{\Gamma(1+g+\hat{\rho}_j+\hat{\rho}_k)\Gamma(1+g+\hat{\rho}_j-\hat{\rho}_k)} \\
&&\makebox[4em]{} \times
\frac{\Gamma(1+\hat{\rho}_j+\hat{\rho}_k)\Gamma(1+\hat{\rho}_j-\hat{\rho}_k)}
{\Gamma(1-g+\hat{\rho}_j+\hat{\rho}_k)\Gamma(1-g+\hat{\rho}_j-\hat{\rho}_k)}
\nonumber \\
&&\times \prod_{1\leq j\leq n}
\frac{\Gamma(1+2\hat{\rho}_j)^2}
     {\prod_{r=1}^{4}
      \Gamma( 1+\hat{g}_r+\hat{\rho}_j) \Gamma( 1-\hat{g}_r+\hat{\rho}_j) }
      \nonumber
\end{eqnarray}
\end{subequations}
for $q=1$,
may be rewritten by canceling common factors in the numerator and
the denominator as
\begin{eqnarray}\label{rhsM}
\lefteqn{\makebox[2em]{}\frac{\hat{\mathcal{C}}_{-,q}(\hat{\rho})}
     {\hat{\mathcal{C}}_{+,q}(\hat{\rho})}= } && \\
&& \!\!\!\!\! \left\{
\begin{array}{ll}
{\displaystyle
\prod_{1\leq j\leq n}
\frac{ (q, q^{1+jg};q)_\infty
       \prod_{1\leq r < s\leq 4} (q^{1+(n-j)g+g_r+g_s};q)_\infty}
     {(q^{1+g},q^{1+(2n-j-1)g+g_1+g_2+g_3+g_4} ;q)_\infty} } &
                                        \text{for}\;\; 0<q<1 \\ [2ex]
{\displaystyle
\prod_{1\leq j\leq n}
\frac{\Gamma (1+g)\,\Gamma(1+(2n-j-1)g+g_1+g_2+g_3+g_4)}
     {\Gamma (1+jg)\, \prod_{1\leq r < s\leq 4} \Gamma (1+(n-j)g+g_r+g_s)} }
& \text{for}\;\; q=1 .
\end{array} \right. \nonumber
\end{eqnarray}
It is clear from this last formula and the explicit expressions
for $\mathcal{C}_{+,q}(\mathbf{x})$ in \eqref{C+q}, \eqref{C+1} that
the summation formula of Theorem~\ref{sumMthm} reduces to
the sums in \eqref{bai-sum}, \eqref{doug-sum} for $n=1$.

\subsection{A generalized Aomoto-Ito sum}
In \cite{aom:product,ito:theta} Aomoto and Ito presented
an evaluation conjecture for certain $q$-Selberg type Jackson integrals
(sums) associated to the reduced root systems.
For the type $A$ root system, the validity of the summation formula
in question was inferred by Aomoto (in \cite{aom:product}) and also
by Kaneko \cite{kan:q-selberg}.
In this special case the sum amounts to an extension
of a $q$-Selberg Jackson integral (sum)
due to Askey, Kadell, Habsieger and Evans
\cite{ask:some,kad:proof,hab:q-integrale,eva:multidimensional}.
For the remaining root systems, the validity of the Aomoto-Ito formula
was verified in the rank two case by Ito, who also
extended the result to the case of a rank two nonreduced $BC$ type
root system \cite{ito:theta}.
Recently, Macdonald observed that the Aomoto-Ito formula follows
for arbitrary reduced root system from what we have dubbed here
`the Macdonald sum' (i.e., the analogue of the sum in \eqref{sumM}
for a reduced root system) \cite{mac:formal}.
The main purpose of the present section is to derive
a generalized Aomoto-Ito type sum for the nonreduced (affine) root
systems. To this end, we will rewrite the sum of Theorem~\ref{sumMthm}
following (in essence) Macdonald's treatment for the case of a reduced root
system.

The main point of our discussion is that multiplication of
both sides of \eqref{sumM} by
the factor $\mathcal{C}_{+,q}(-\mathbf{z})/\mathcal{C}_{-,q}(\mathbf{z})$
leads to

{\em A generalized Aomoto-Ito sum}
\begin{equation}\label{sumAI}
\sum_{\lambda \in \mathbb{Z}^n}
\frac{1}
     {\mathcal{C}_{+,q}(\mathbf{z}+\lambda )\,
      \mathcal{C}_{-,q}(\mathbf{z}+\lambda ) }
=
\frac{\mathcal{C}_{+,q}(-\mathbf{z})}
     {\mathcal{C}_{-,q}(\mathbf{z}) }
\frac{\hat{\mathcal{C}}_{-,q}(\hat{\rho})}
     {\hat{\mathcal{C}}_{+,q}(\hat{\rho})} ,
\end{equation}
where it
is again assumed that the parameters $g$, $g_r$ ($r=1,2,3,4$)
satisfy the convergence condition in \eqref{parco}.
To avoid poles one should choose $\mathbf{z}$ such that
$-g+z_j\pm z_k$ ($1\leq j<k\leq n$) and $-g_r+z_j$
($r=1,2,3,4$; $1\leq j\leq n$) are nonzero modulo
the lattice $\Omega_q$ \eqref{lattice} (cf. the expressions below).
In fact, the identity \eqref{sumAI} may be viewed as an equality
between meromorphic functions of $\mathbf{z}$. More specifically,
the l.h.s. and r.h.s. of \eqref{sumAI} are equal as a
meromorphic function in $z_j$ (where $j$ is arbitrary but fixed)
with poles congruent (modulo $\Omega_q$ \eqref{lattice})
to $z_j=\pm (g-z_k)$ ($k<j$), $z_j=g\pm z_k$
($k>j$) and $z_j=g_r$ ($r=1,2,3,4$).

In order to arrive at the Aomoto-Ito type sum \eqref{sumAI} from \eqref{sumM},
we have used for the l.h.s. the fact that the multiplier
$\mathcal{C}_{+,q}(-\mathbf{z})/\mathcal{C}_{-,q}(\mathbf{z})$
is periodic in $\mathbf{z}$
\begin{equation}\label{periodicity}
\frac{\mathcal{C}_{+,q}(-\mathbf{z} )}
     {\mathcal{C}_{-,q}(\mathbf{z} ) }=
\frac{\mathcal{C}_{+,q}(-\mathbf{z}-\lambda )}
     {\mathcal{C}_{-,q}(\mathbf{z}+\lambda ) },
\;\;\;\;\;\;\;\;\; \text{for}\;\;\; \lambda \in \mathbb{Z}^n.
\end{equation}
This periodicity is not difficult to see by inferring that the factor
$\mathcal{C}_{+,q}(-\mathbf{z})/\mathcal{C}_{-,q}(\mathbf{z})$ is invariant
with respect to translations in $\mathbf{z}$ over the
unit vectors $e_j$, $j=1,\ldots ,n$
(of the standard basis in $\mathbb{R}^n$)
using the elementary shift properties
for the $q$-shifted factorial (viz. $(qa;q)_\infty = (a;q)_\infty / (1-a)$)
and gamma function (viz. $\Gamma (a+1)=a\,\Gamma(a)$), respectively.
(For a still simpler way to deduce the periodicity see below.)
For $0<q<1$ and special values of the parameters $g_1,g_2,g_3,g_4$,
the summation formula in \eqref{sumAI} amounts to
the $BC$ type Aomoto-Ito sum (see the appendix of \cite{ito:theta}).
By further specialization
of the parameters one recovers the Aomoto-Ito sums associated
to the reduced root systems of type $B$, $C$ and $D$.

To make the contact with \cite{aom:product,ito:theta} and
\cite{mac:formal} more explicit, it is helpful to observe
that for $0<q<1$ the terms in \eqref{sumAI} may be written as
\begin{eqnarray}\label{aomtermq}
\lefteqn{\frac{1}
     {\mathcal{C}_{+,q}(\mathbf{x})\,
      \mathcal{C}_{-,q}(\mathbf{x}) } =
q^{\sum_{j= 1}^n (1+2\hat{\rho}_j) x_j} } &&  \\
&&  \times
\prod_{1\leq j<k\leq n}
(1-q^{x_j+x_k})(1-q^{x_j-x_k})
\frac{(q^{1+g+x_j+x_k},q^{1+g+x_j-x_k};q)_\infty}
     {(q^{-g+x_j+x_k},q^{-g+x_j-x_k};q)_\infty} \nonumber \\
&& \times \prod_{1\leq j\leq n}
(1-q^{2x_j})
\frac{(q^{1+g_1+x_j},q^{1+g_2+x_j},q^{1+g_3+x_j},q^{1+g_4+x_j};q)_\infty}
     {(q^{-g_1+x_j},q^{-g_2+x_j},q^{-g_3+x_j},q^{-g_4+x_j};q)_\infty}
\nonumber
\end{eqnarray}
(with $\mathbf{x}=\mathbf{z}+\lambda$).
Furthermore, the periodic Aomoto factor
$\mathcal{C}_{+,q}(-\mathbf{z})/\mathcal{C}_{-,q}(\mathbf{z})$
on the r.h.s. is conveniently
rewritten in terms of the Jacobi theta function $\theta$ \eqref{theta}
as follows
\begin{eqnarray}\label{aomF}
\lefteqn{\makebox[2em]{}
\frac{\mathcal{C}_{+,q}(-\mathbf{z})}{\mathcal{C}_{-,q}(\mathbf{z})}
=q^{\sum_{j= 1}^n (1+2\hat{\rho}_j) z_j} } && \\
&& \times
\prod_{1\leq j<k\leq n}
\frac{\theta (q^{z_j+z_k})\, \theta (q^{z_j-z_k})}
     {\theta (q^{-g+z_j+z_k})\,\theta(q^{-g+z_j-z_k})}
\prod_{1\leq j\leq n}
\frac{(q;q)_\infty^3\theta (q^{2z_j})}
     {\prod_{1\leq r\leq 4} \theta (q^{-g_r+z_j})} \nonumber
\end{eqnarray}
(this is clear from the product representation for $\theta (\zeta)$
in \eqref{prod}). Observe that the periodicity relation
\eqref{periodicity} is now easily deduced using the
quasi-periodicity relation \eqref{quasi} for the theta function.

For $q=1$ the summation formula in \eqref{sumAI} constitutes
a degeneration of the (generalized) Aomoto-Ito formula.
The summand on the l.h.s. is in this case governed by
\begin{eqnarray}\label{aomterm1}
\lefteqn{\makebox[2em]{}\frac{1}
     {\mathcal{C}_{+,1}(\mathbf{x})\,
      \mathcal{C}_{-,1}(\mathbf{x}) }
=} && \\
&& \prod_{1\leq j<k\leq n}
\frac{(x_j+x_k)(x_j-x_k)\,\Gamma (-g+x_j+x_k)\, \Gamma (-g+x_j-x_k)}
     {\Gamma (1+g+x_j+x_k)\, \Gamma (1+g+x_j-x_k)} \nonumber \\
&& \times \prod_{1\leq j\leq n}
\frac{2x_j\, \Gamma (-g_1+x_j)\, \Gamma (-g_3+x_j)\,
             \Gamma (-g_3+x_j)\,\Gamma (-g_4+x_j) }
     {\Gamma (1+g_1+x_j)\,\Gamma (1+g_2+x_j)\,
      \Gamma (1+g_3+x_j)\, \Gamma (1+g_4+x_j)} \nonumber
\end{eqnarray}
($\mathbf{x}=\mathbf{z}+\lambda$)
and the degenerate Aomoto factor on the r.h.s. can now be
written in the manifestly periodic form
\begin{eqnarray}
\frac{\mathcal{C}_{+,1}(-\mathbf{z})}{\mathcal{C}_{-,1}(\mathbf{z})}
 &=&
\prod_{1\leq j<k\leq n}
\frac{\sin\pi (z_j+z_k)\,\sin\pi (z_j-z_k)}
     {\sin\pi (-g+z_j+z_k)\,\sin\pi (-g+z_j-z_k)} \\
&& \times \prod_{1\leq j\leq n}
\frac{\pi^3\sin (2\pi z_j)}
     {\prod_{1\leq r\leq 4} \sin \pi (-g_r+z_j)} \nonumber
\end{eqnarray}
(where we have used the reflection equation \eqref{refg}).

For $n=1$ the generalized Aomoto-Ito
summation formula
reduces (after division by a common factor $q^{(1+g_1+g_2+g_3+g_4)z}$) to
\begin{eqnarray}
\lefteqn{\sum_{\lambda\in\mathbb{Z}}
q^{(1+g_1+g_2+g_3+g_4)\lambda}
(1-q^{2z+2\lambda})
\prod_{1\leq r\leq 4} \frac{(q^{1+g_r+z+\lambda};q)_\infty}
                           {(q^{-g_r+z+\lambda};q)_\infty} } && \\
&& \makebox[4em]{}=
\frac{(q;q)_\infty^3\theta (q^{2z})}
     {\prod_{1\leq r\leq 4} \theta (q^{-g_r+z})}
\frac{ (q ;q)_\infty \prod_{1\leq r < s\leq 4} (q^{1+g_r+g_s};q)_\infty}
     {(q^{1+g_1+g_2+g_3+g_4} ;q)_\infty}. \nonumber
\end{eqnarray}
and the corresponding $q=1$ degeneration reads
\begin{eqnarray}
\lefteqn{
\sum_{\lambda\in\mathbb{Z}}
\Bigl(
(2z+2\lambda )\prod_{1\leq r\leq 4}
\frac{\Gamma (-g_r+z+\lambda )}
     {\Gamma (1+g_3+z+\lambda )}  \Bigr) } &&   \\
&& \makebox[5em]{} =
\frac{\pi^3\sin (2\pi z_j)}
     {\prod_{1\leq r\leq 4} \sin \pi (-g_r+z_j)}
\frac{\Gamma(1+g_1+g_2+g_3+g_4)}
     {\prod_{1\leq r < s\leq 4} \Gamma (1+g_r+g_s)} . \nonumber
\end{eqnarray}

\subsection{A multiple Bailey ${}_6\psi_6$ and
Dougall ${}_5H_5$ sum}\label{sec2.4}
Division of the Macdonald type sum \eqref{sumM} by the
middle term
$1/(\mathcal{C}_{+,q}(\mathbf{z} )\,\mathcal{C}_{+,q}(-\mathbf{z}) )$
(corresponding to $\lambda =\mathbf{0}$),
or equivalently, division of the Aomoto-Ito type sum \eqref{sumAI}
by the middle term
$1/(\mathcal{C}_{+,q}(\mathbf{z} )\,\mathcal{C}_{-,q}(\mathbf{z}) )$
leads us to

{\em A multiple Bailey/Dougall sum}
\begin{subequations}
\begin{equation}\label{sumBDa}
\sum_{\lambda \in \mathbb{Z}^n}
\frac{\mathcal{C}_{+,q}(\mathbf{z} )\,
      \mathcal{C}_{+,q}(-\mathbf{z} )}
     {\mathcal{C}_{+,q}(\mathbf{z}+\lambda )\,
      \mathcal{C}_{+,q}(-\mathbf{z}-\lambda ) }
=
\mathcal{C}_{+,q}(\mathbf{z} )\,
      \mathcal{C}_{+,q}(-\mathbf{z} )\:
\frac{\hat{\mathcal{C}}_{-,q}(\hat{\rho})}
     {\hat{\mathcal{C}}_{+,q}(\hat{\rho})}
\end{equation}
or equivalently
\begin{equation}\label{sumBDb}
\sum_{\lambda \in \mathbb{Z}^n}
\frac{\mathcal{C}_{+,q}(\mathbf{z} )\,
      \mathcal{C}_{-,q}(\mathbf{z} )}
     {\mathcal{C}_{+,q}(\mathbf{z}+\lambda )\,
      \mathcal{C}_{-,q}(\mathbf{z}+\lambda ) }
=
\mathcal{C}_{+,q}(\mathbf{z} )\,
      \mathcal{C}_{+,q}(-\mathbf{z} )\:
\frac{\hat{\mathcal{C}}_{-,q}(\hat{\rho})}
     {\hat{\mathcal{C}}_{+,q}(\hat{\rho})}
\end{equation}
\end{subequations}
(the equality of the terms on the l.h.s. of \eqref{sumBDa}
and \eqref{sumBDb} traces back to the periodicity relation
\eqref{periodicity}), where the parameters
$g$, $g_r$ again satisfy the convergence condition in \eqref{parco}.
To avoid poles, the components of $\mathbf{z}$ should now be
chosen subject to the genericity condition
that $z_j$, $z_j\pm z_k$ ($j\neq k$) is nonzero
(modulo $\frac{2\pi}{i\log q}$)
and that $g_r\pm z_j$, $g\pm z_j\pm z_k$ ($j\neq k$)
is not a negative integer (modulo $\frac{2\pi}{i\log q}$).
In order to facilitate the comparison
with the classical one-variable Bailey and Dougall sums
in \eqref{bai-ser} and \eqref{doug-ser}, it is useful
to display the terms on the l.h.s. of the multiple summation
formula \eqref{sumBDa}, \eqref{sumBDb} somewhat more explicitly.
One has for $0<q<1$ that
\begin{subequations}
\begin{eqnarray}\label{baiterm}
\lefteqn{\makebox[2em]{}
\frac{\mathcal{C}_{+,q}(\mathbf{z} )\,
      \mathcal{C}_{-,q}(\mathbf{z} )}
     {\mathcal{C}_{+,q}(\mathbf{z}+\lambda )\,
      \mathcal{C}_{-,q}(\mathbf{z}+\lambda ) } =
     q^{\sum_{j=1}^n (1+2\hat{\rho}_j) \lambda_j}
                                              } && \\
&& \!\!\!\!\! \times \prod_{1\leq j<k\leq n}
\Bigl( \frac{1-q^{z_j+z_k+\lambda_j+\lambda_k}}
            {1-q^{z_j+z_k}}
      \frac{1-q^{z_j-z_k+\lambda_j-\lambda_k}}
            {1-q^{z_j-z_k}} \Bigr)  \nonumber \\
&& \makebox[5em]{}\times \frac{(q^{-g+z_j+z_k};q)_{\lambda_j+\lambda_k}\,
      (q^{-g+z_j-z_k};q)_{\lambda_j-\lambda_k} }
     {(q^{1+g+z_j+z_k};q)_{\lambda_j+\lambda_k}\,
      (q^{1+g+z_j-z_k};q)_{\lambda_j-\lambda_k} }    \nonumber \\
&& \!\!\!\!\! \times \prod_{1\leq j\leq n}
\Big( \frac{1-q^{2z_j+2\lambda_j}}
           {1-q^{2z_j}} \Bigr)
\frac{(q^{-g_1+z_j},q^{-g_2+z_j},q^{-g_3+z_j},q^{-g_4+z_j};q)_{\lambda_j}}
   {(q^{1+g_1+z_j},q^{1+g_2+z_j},q^{1+g_3+z_j},q^{1+g_4+z_j};q)_{\lambda_j}}
     \nonumber
\end{eqnarray}
(using $(a;q)_\infty /(q^m a;q)_\infty = (a;q)_m$) and for $q=1$ that
\begin{eqnarray}\label{dougterm}
\lefteqn{\makebox[1em]{}
\frac{\mathcal{C}_{+,1}(\mathbf{z} )\,
      \mathcal{C}_{-,1}(\mathbf{z} )}
     {\mathcal{C}_{+,1}(\mathbf{z}+\lambda )\,
      \mathcal{C}_{-,1}(\mathbf{z}+\lambda ) } =
                                              } && \\
&& \prod_{1\leq j< k\leq n}
\Bigl( 1+ \frac{\lambda_j+\lambda_k}{z_j+z_k} \Bigr)
\Bigl( 1+ \frac{\lambda_j-\lambda_k}{z_j-z_k} \Bigr)\nonumber \\
&& \makebox[3em]{}
\times \frac{(-g+z_j+z_k)_{\lambda_j+\lambda_k} \,
             (-g+z_j-z_k)_{\lambda_j-\lambda_k}}
            {(1+g+z_j+z_k)_{\lambda_j+\lambda_k} \,
             (1+g+z_j-z_k)_{\lambda_j-\lambda_k}}
\nonumber \\
&&\hspace{-4em}\times \prod_{1\leq j\leq n}
\Bigl( 1+\frac{\lambda_j}{z_j}\Bigr)
\frac{(-g_1+z_j,-g_2+z_j,-g_3+z_j,-g_4+z_j)_{\lambda_j}}
     {(1+g_1+z_j,1+g_2+z_j,1+g_3+z_j,1+g_4+z_j)_{\lambda_j}}  \nonumber
\end{eqnarray}
\end{subequations}
(using $\Gamma (a+m)/\Gamma (a) =(a)_m$).
{}From \eqref{baiterm}, \eqref{dougterm} and \eqref{rhsM} it is immediate
that the summation formula \eqref{sumBDa}/\eqref{sumBDb}
specializes for $n=1$ to the Bailey and Dougall sums in
\eqref{bai-ser} and \eqref{doug-ser}.

\begin{remark}
Multidimensional analogues of the
Bailey ${}_6\psi_6$ and Dougall ${}_5F_4$
summation formulas \eqref{bai-ser}, \eqref{doug-ser}
different from the ones in \eqref{sumBDa}/\eqref{sumBDb}
were introduced by Gustafson in \cite{gus:multilateral,gus:macdonald}
(cf. also Section~\ref{sec4}) and
recently still another multidimensional
version (in two distinct variations) of the ${}_6\psi_6$ sum
was presented by Schlosser \cite{sch:summation}.
For $g=-1/2$ (and $0<q<1$),
the sum in \eqref{sumBDa}/\eqref{sumBDb} may be seen as a special case of
Schlosser's multiple Bailey sums.
\end{remark}

\section{Truncation: multiple analogues of very-well-poised unilateral
(basic) hypergeometric series}\label{sec3}

In this section the vector $\mathbf{z}$ will be specialized
in such a way that the multiple Baily/Dougall sum of Section~\ref{sec2.4}
truncates to a sum over the dominant cone
\begin{equation}\label{cone}
\Lambda = \{ \lambda\in \mathbb{Z}^n\; |\;
           \lambda_1\geq \lambda_2\geq \cdots \geq \lambda_n\geq 0 \; \}
\end{equation}
or---for special choice of the parameters---to
a finite sum over the sub-alcove
\begin{equation}\label{alcove}
\Lambda_N = \{ \lambda\in \mathbb{Z}^n\; |\; N\geq
           \lambda_1\geq \lambda_2\geq \cdots \geq \lambda_n\geq 0 \; \}
\;\;\;\;\;\;\;\;\; (\text{with}\;\; N\in\mathbb{N}).
\end{equation}
The resulting summation formulas constitute multiple analogues of
the (terminating) very-well-poised Rogers ${}_6\phi_5$ and Dougall ${}_5F_4$
sums \cite{rog:third,dou:vandermonde,gas-rah:basic}.
For $0<q<1$ the terms of the multiple series in question are of the form
\begin{subequations}
\begin{eqnarray}\label{rogterm}
\makebox[1em]{}
\Delta_q (\lambda ) \!\!\!\! &=& \!\!\!
     q^{\sum_{j=1}^n (1-2\hat{\rho}_j) \lambda_j}  \\
&& \!\!\!\!\! \times \prod_{1\leq j<k\leq n}
\Bigl( \frac{1-q^{\rho_j+\rho_k+\lambda_j+\lambda_k}}
            {1-q^{\rho_j+\rho_k}}
      \frac{1-q^{\rho_j-\rho_k+\lambda_j-\lambda_k}}
            {1-q^{\rho_j-\rho_k}} \Bigr)  \nonumber \\
&& \makebox[5em]{}\times \frac{(q^{g+\rho_j+\rho_k};q)_{\lambda_j+\lambda_k}\,
      (q^{g+\rho_j-\rho_k};q)_{\lambda_j-\lambda_k} }
     {(q^{1-g+\rho_j+\rho_k};q)_{\lambda_j+\lambda_k}\,
      (q^{1-g+\rho_j-\rho_k};q)_{\lambda_j-\lambda_k} }    \nonumber \\
&& \!\!\!\!\! \times \prod_{1\leq j\leq n}
\Big( \frac{1-q^{2\rho_j+2\lambda_j}}
           {1-q^{2\rho_j}} \Bigr)
\frac{\prod_{r= 1}^4(q^{g_r+\rho_j};q)_{\lambda_j}}
   {\prod_{r= 1}^4(q^{1-g_r+\rho_j};q)_{\lambda_j}}
     \nonumber
\end{eqnarray}
(corresponding to a Rogers type series) and for $q=1$ we have
\begin{eqnarray}\label{douterm}
\makebox[1em]{} \Delta_1 (\lambda )\!\!\!\! &=&\!\!\!\!\!\!
\prod_{1\leq j< k\leq n}
\Bigl( 1+ \frac{\lambda_j+\lambda_k}{\rho_j+\rho_k} \Bigr)
\Bigl( 1+ \frac{\lambda_j-\lambda_k}{\rho_j-\rho_k} \Bigr)  \\
&& \makebox[3em]{}
\times \frac{(g+\rho_j+\rho_k)_{\lambda_j+\lambda_k} \,
             (g+\rho_j-\rho_k)_{\lambda_j-\lambda_k}}
            {(1-g+\rho_j+\rho_k)_{\lambda_j+\lambda_k} \,
             (1-g+\rho_j-\rho_k)_{\lambda_j-\lambda_k}}
\nonumber \\
&&\times \prod_{1\leq j\leq n}
\Bigl( 1+\frac{\lambda_j}{\rho_j}\Bigr)
\frac{\prod_{r= 1}^4 (g_r+\rho_j)_{\lambda_j}}
     {\prod_{r= 1}^4 (1-g_r+\rho_j)_{\lambda_j}}
\nonumber
\end{eqnarray}
\end{subequations}
(corresponding to a Dougall type series),
where $\rho_j$ and $\hat{\rho}_j$ are given by \eqref{rho}.

\subsection{Multiple Rogers ${}_6\phi_5$ and Dougall ${}_5F_4$ sums}
It is not so difficult to infer that the terms of the Bailey/Dougall
sum \eqref{sumBDa}, \eqref{sumBDb} become zero for
$\lambda \in \mathbb{Z}^n\setminus \Lambda$ (with $\Lambda$ given
by \eqref{cone}) if one sets $\mathbf{z}=-\rho$.
Indeed, when $z_j=-\rho_j$ \eqref{rho}
($j=1,\ldots ,n$) the factors
$(q^{1+g+z_j-z_{j+1}};q)_{\lambda_j-\lambda_{j+1}}$,
$(1+g+z_j-z_{j+1})_{\lambda_j-\lambda_{j+1}}$ and
$(q^{1+g_a+z_n};q)_{\lambda_n}$,
$(1+g_a+z_n)_{\lambda_n}$
in the denominators
of \eqref{baiterm}, \eqref{dougterm} give rise to a zero
for $\lambda_j <\lambda_{j+1}$ and $\lambda_n<0$, respectively.
Simplification of the r.h.s. of \eqref{sumBDa}/\eqref{sumBDb}
(for $\mathbf{z}=-\rho$ the factors $\mathcal{C}_{+,q}(\rho )$
and $1/\hat{\mathcal{C}}_{+,q}(\hat{\rho})$ cancel each other)
and reflection of the parameters in the origin ($g\rightarrow -g$,
$g_r\rightarrow -g_r$) so as to avoid the excessive appearance
of minus signs, finally leads us to the following summation theorem.

\begin{theorem}\label{sumRDthm}
Let $0<q\leq 1$.
For parameters subject to the condition
\begin{equation}\label{convcond}
\text{Re}\: (1-2(n-j)g-g_1-g_2-g_3-g_4) >0
\end{equation}
($j=1,\ldots ,n$),
one has that
\begin{equation}\label{sumRD}
\sum_{\lambda\in\Lambda} \Delta_q (\lambda )
= \mathcal{N}_{q,\Lambda}
\end{equation}
with
\begin{subequations}
\begin{eqnarray}\label{rognorm}
\makebox[1em]{}\mathcal{N}_{q,\Lambda}
\!\!\!\!\! &=&\!\!\!\!\!
\prod_{1\leq j<k\leq n}
\frac{(q^{1+\rho_j+\rho_k},
       q^{1+\rho_j-\rho_k},
       q^{1+g-\hat{\rho}_j-\hat{\rho}_k},
       q^{1+g-\hat{\rho}_j+\hat{\rho}_k};q)_\infty}
     {(q^{1-g+\rho_j+\rho_k},
       q^{1-g+\rho_j-\rho_k},
       q^{1-\hat{\rho}_j-\hat{\rho}_k},
       q^{1-\hat{\rho}_j+\hat{\rho}_k};q)_\infty} \\
&&\times \prod_{1\leq j\leq n}
\frac{(q^{1+2\rho_j};q)_\infty
       \prod_{r=1}^{4} ( q^{1+\hat{g}_r-\hat{\rho}_j};q)_\infty}
     {(q^{1-2\hat{\rho}_j};q)_\infty
       \prod_{r=1}^{4} ( q^{1-g_r+\rho_j};q)_\infty}  \nonumber
\end{eqnarray}
for $0<q<1$ and
\begin{eqnarray}\label{dougnorm}
\makebox[1em]{}\mathcal{N}_{1,\Lambda}
\!\!\!\!\!&=&\!\!\!\!\!
\prod_{1\leq j<k\leq n}
\frac{\Gamma(1-g+\rho_j+\rho_k)\Gamma(1-g+\rho_j-\rho_k)}
     {\Gamma(1+\rho_j+\rho_k)\Gamma(1+\rho_j-\rho_k)} \\
&&\makebox[4em]{} \times
\frac{\Gamma(1-\hat{\rho}_j-\hat{\rho}_k)\Gamma(1-\hat{\rho}_j+\hat{\rho}_k)}
{\Gamma(1+g-\hat{\rho}_j-\hat{\rho}_k)\Gamma(1+g-\hat{\rho}_j+\hat{\rho}_k)}
\nonumber \\
&&\times \prod_{1\leq j\leq n}
\frac{\Gamma(1-2\hat{\rho}_j)
      \prod_{r=1}^{4} \Gamma( 1-g_r+\rho_j)}
     {\Gamma (1+2\rho_j)
      \prod_{r=1}^{4} \Gamma( 1+\hat{g}_r-\hat{\rho}_j)}.  \nonumber
\end{eqnarray}
\end{subequations}
Moreover, the series on the l.h.s. of \eqref{sumRD}
converges in absolute value.
\end{theorem}

By canceling common terms in numerator and denominator the
evaluation constants \eqref{rognorm} and \eqref{dougnorm}
may be rewritten as
\begin{subequations}
\begin{eqnarray}
\mathcal{N}_{q,\Lambda}
\!\!\!\!\! &=&\!\!\!\!\!
\prod_{1\leq j\leq n}
\frac{(q^{1+(2n-j-1)g+2g_a};q)_\infty}
     {(q^{1-(2n-j-1)g-g_a-g_b-g_c-g_d};q)_\infty} \\
&& \makebox[4em]{}\times
\frac{ \prod_{\stackrel{1\leq r<s\leq 4}{r,s\neq a}}
(q^{1-(n-j)g-g_r-g_s};q)_\infty }
     { \prod_{\stackrel{1\leq r\leq 4}{r\neq a}}
(q^{1+(n-j)g+g_a-g_r};q)_\infty }
\nonumber
\end{eqnarray}
($0<q<1$) and
\begin{eqnarray}
\mathcal{N}_{1,\Lambda}
\!\!\!\!\! &=&\!\!\!\!\!
\prod_{1\leq j\leq n}
\frac{\Gamma (1-(2n-j-1)g-g_a-g_b-g_c-g_d)}{\Gamma (1+(2n-j-1)g+2g_a)} \\
&& \makebox[4em]{}\times
\frac{\prod_{\stackrel{1\leq r\leq 4}{r\neq a}} \Gamma (1+(n-j)g+g_a-g_r)}
     {\prod_{\stackrel{1\leq r<s\leq 4}{r,s\neq a}} \Gamma (1-(n-j)g-g_r-g_s)}
,
\nonumber
\end{eqnarray}
\end{subequations}
respectively.
When $n=1$ the summation formula in Theorem~\ref{sumRDthm}
reduces for $0<q<1$ to the nonterminating
Rogers sum (cf. \cite{rog:third,gas-rah:basic})
\begin{subequations}
\begin{eqnarray}\label{rog-uni}
\lefteqn{
\makebox[2em]{}\sum_{\lambda \in \mathbb{N}}
q^{(1-g_a-g_b-g_c-g_d)\lambda}
\Big( \frac{1-q^{2g_a+2\lambda}}
           {1-q^{2g_a}} \Bigr)
\prod_{1\leq r\leq 4}
\frac{(q^{g_r+g_a};q)_\lambda}
     {(q^{1-g_r+g_a};q)_\lambda} } &&  \\
&&={}_6\phi_5 \Bigl(
\begin{array}{c}
q^{1+g_a},-q^{1+g_a}, q^{2g_a},q^{g_a+g_b},q^{g_a+g_c},q^{g_a+g_d} \\
q^{g_a}, -q^{g_a}, q^{1+g_a-g_b},q^{1+g_a-g_c},q^{1+g_a-g_d}
\end{array} ;q, q^{1-g_a-g_b-g_c-g_d} \Bigr)  \nonumber \\
&& =
\frac{ (q^{1+2g_a}, q^{1-g_b-g_c},q^{1-g_b-g_d},q^{1-g_c-g_d};q)_\infty}
 {(q^{1+g_a-g_b},q^{1+g_a-g_c},q^{1+g_a-g_d},q^{1-g_a-g_b-g_c-g_d} ;q)_\infty}
\nonumber
\end{eqnarray}
and for $q=1$ to the nonterminating Dougall sum (cf.
\cite{dou:vandermonde,gas-rah:basic})
\begin{eqnarray}\label{doug-uni}
\lefteqn{\makebox[1em]{}
\sum_{\lambda \in \mathbb{N}}\:
 \Bigl( 1+\frac{\lambda}{g_a}\Bigr)
\prod_{1\leq r\leq 4}\frac{(g_r+g_a)_\lambda}
                          {(1-g_r+g_a)_\lambda} } && \\
&& \hspace{-2em}=
{}_5F_4 \left(
\begin{array}{c}
1+g_a, 2g_a,g_a+g_b,g_a+g_c,g_a+g_d \\
g_a,1+g_a-g_b,1+g_a-g_c,1+g_a-g_d
\end{array}; 1 \right) \nonumber \\
&&\hspace{-2em} =
\frac{\Gamma (1+g_a-g_b)\Gamma (1+g_a-g_c)\Gamma (1+g_a-g_d)
                                   \Gamma(1-g_a-g_b-g_c-g_d)}
     {\Gamma (1+2g_a)\Gamma (1-g_b-g_c)\Gamma (1-g_b-g_d)\Gamma (1-g_c-g_d)} ,
 \nonumber
\end{eqnarray}
\end{subequations}
where the parameters are assumed
to satisfy the convergence condition $\text{Re}\: (1-g_a-g_b-g_c-g_d)>0$.

\subsection{Terminating multiple Rogers ${}_6\phi_5$ and
Dougall ${}_5F_4$ sums}
When the parameters in the Rogers/Dougall type sum of
Theorem~\ref{sumRDthm} are chosen in such
a way that $(n-1)g+g_a+g_b+N=0$
with $N\in\mathbb{N}$, then the series on the l.h.s. terminates as the terms
$\Delta_q(\lambda)$ \eqref{rogterm}, \eqref{douterm}
become zero for $\lambda\in \Lambda\setminus\Lambda_N$
(where $\Lambda_N$ is given by \eqref{alcove}).
This is because we now
pick up a zero from the factor $(q^{g_b+\rho_1};q)_{\lambda_1}$ or
$(g_b+\rho_1)_{\lambda_1}$ in the numerator when $\lambda_1>N$.
As a consequence, the l.h.s. of \eqref{sumRD} becomes a rational
expression in $q$, $q^g$, $q^{g_r}$ (for $q\neq 1$) or
in $g$, $g_r$ (for $q=1$) and the same must be true for the corresponding
r.h.s. Indeed, for the parameters subject to the above truncation
condition the infinite products entering
$\mathcal{N}_{\Lambda}$ \eqref{rognorm}, \eqref{dougnorm}
may be reduced to finite
products by canceling common factors in the numerator and denominator.
This way one arrives at an expression for the r.h.s. that can be written as
\begin{subequations}
\begin{eqnarray}\label{finrognorm}
\makebox[2em]{}\mathcal{N}_{q,\Lambda_N}
\!\!\!\! &=&\!\!\!\!\!\!\!
\prod_{1\leq j<k\leq n}
\frac{(q^{1+\rho_j+\rho_k};q)_N }
     {(q^{1-g+\rho_j+\rho_k};q)_N }
\prod_{1\leq j\leq n}
\frac{(q^{1+2\rho_j},q^{1+\hat{g}_b-\hat{\rho}_j};q)_N}
     {(q^{1-g_c+\rho_j},q^{1-g_d+\rho_j};q)_N} \\
\!\!\!\! &=&\!\!\!\!\!\!\!
\prod_{1\leq j<k\leq n}
\frac{(q^{1+g-\hat{\rho}_j-\hat{\rho}_k};q)_{-N} }
     {(q^{1-\hat{\rho}_j-\hat{\rho}_k};q)_{-N} }
\prod_{1\leq j\leq n}
\frac{(q^{1+\hat{g}_c-\hat{\rho}_j},q^{1+\hat{g}_d-\hat{\rho}_j};q)_{-N}}
     {(q^{1-2\hat{\rho}_j},q^{1-g_b+\rho_j};q)_{-N}}
\nonumber
\end{eqnarray}
for $q\neq 1$ and as
\begin{eqnarray}\label{findougnorm}
\makebox[2em]{}\mathcal{N}_{1,\Lambda_N}
\!\!\!\! &=&\!\!\!\!\!\!\!
\prod_{1\leq j<k\leq n}
\frac{(1+\rho_j+\rho_k)_N }
     {(1-g+\rho_j+\rho_k)_N }
\prod_{1\leq j\leq n}
\frac{(1+2\rho_j,1+\hat{g}_b-\hat{\rho}_j)_N}
     {(1-g_c+\rho_j,1-g_d+\rho_j)_N}                  \\
\!\!\!\! &=&\!\!\!\!\!\!\!
\prod_{1\leq j<k\leq n}
\frac{(1+g-\hat{\rho}_j-\hat{\rho}_k)_{-N} }
     {(1-\hat{\rho}_j-\hat{\rho}_k)_{-N} }
\prod_{1\leq j\leq n}
\frac{(1+\hat{g}_c-\hat{\rho}_j,1+\hat{g}_d-\hat{\rho}_j)_{-N}}
     {(1-2\hat{\rho}_j,1-g_b+\rho_j)_{-N}}
\nonumber
\end{eqnarray}
\end{subequations}
for $q=1$.
For the equality of the expressions on the first line and the second
line of the formulas \eqref{finrognorm}, \eqref{findougnorm} it is essential
that
the parameters satisfy the truncation condition $(n-1)g+g_a+g_b+N=0$.
The corresponding expression for the r.h.s. of the terminating
multiple Rogers/Dougall summation formula becomes somewhat more
symmetric if we combine the both representations for
$\mathcal{N}_{q,\Lambda_N}$ into a single expression for
$\mathcal{N}_{q,\Lambda_N}^2$.

\begin{theorem}\label{finsumRDthm}
Let $N\in\mathbb{N}$.
For parameters subject to the truncation condition
\begin{equation}\label{trunc}
(n-1)g+g_a+g_b+ N =0,
\end{equation}
one has that
\begin{equation}\label{finsumRD}
\sum_{\lambda\in\Lambda_N} \Delta_q(\lambda ) =
\mathcal{N}_{q,\Lambda_N}
\end{equation}
with
\begin{subequations}
\begin{eqnarray}
\makebox[1em]{}\mathcal{N}_{q,\Lambda_N}^2
\!\!\!\! &=&\!\!\!\!\!\!\!
\prod_{1\leq j<k\leq n}
\frac{(q^{1+\rho_j+\rho_k};q)_N }
     {(q^{1-g+\rho_j+\rho_k};q)_N }
\frac{(q^{1+g-\hat{\rho}_j-\hat{\rho}_k};q)_{-N} }
     {(q^{1-\hat{\rho}_j-\hat{\rho}_k};q)_{-N} }
\\
&& \hspace{-4em}\times
\prod_{1\leq j\leq n}
\frac{(q^{1+2\rho_j},q^{1+\hat{g}_b-\hat{\rho}_j};q)_N}
     {(q^{1-2\hat{\rho}_j},q^{1-g_b+\rho_j};q)_{-N}}
\frac{(q^{1+\hat{g}_c-\hat{\rho}_j},q^{1+\hat{g}_d-\hat{\rho}_j};q)_{-N}}
     {(q^{1-g_c+\rho_j},q^{1-g_d+\rho_j};q)_N}
\nonumber
\end{eqnarray}
for $q\neq 1$ and
\begin{eqnarray}
\makebox[1em]{}\mathcal{N}_{1,\Lambda_N}^2
\!\!\!\! &=&\!\!\!\!\!\!\!
\prod_{1\leq j<k\leq n}
\frac{(1+\rho_j+\rho_k)_N }
     {(1-g+\rho_j+\rho_k)_N }
\frac{(1+g-\hat{\rho}_j-\hat{\rho}_k)_{-N} }
     {(1-\hat{\rho}_j-\hat{\rho}_k)_{-N} } \\
&& \hspace{-4em}\times
\prod_{1\leq j\leq n}
\frac{(1+2\rho_j,1+\hat{g}_b-\hat{\rho}_j)_N}
     {(1-2\hat{\rho}_j,1-g_b+\rho_j)_{-N}}
\frac{(1+\hat{g}_c-\hat{\rho}_j,1+\hat{g}_d-\hat{\rho}_j)_{-N}}
     {(1-g_c+\rho_j,1-g_d+\rho_j)_N} .
\nonumber
\end{eqnarray}
\end{subequations}
The summation formula \eqref{finsumRD} holds as a rational identity
in $q^g$, $q^{g_a},q^{g_b},q^{g_c},q^{g_d}$ and $q$ (for $q\neq 1$)
or in $g$, $g_a,g_b,g_c,g_d$ (for $q=1$),
subject to the relation \eqref{trunc}.
\end{theorem}

By further cancellation of common factors in the numerator and denominator
the evaluation constants \eqref{finrognorm}, \eqref{findougnorm}
may be rewritten as
\begin{subequations}
\begin{eqnarray}
\makebox[1em]{}\mathcal{N}_{q,\Lambda_N}
\!\!\!\! &=&\!\!\!\!\!\!\!
\prod_{1\leq j\leq n}
\frac{(q^{1+(2n-j-1)g+2g_a},q^{1-(n-j)g-g_c-g_d};q)_N}
     {(q^{1+(n-j)g+g_a-g_c},q^{1+(n-j)g+g_a-g_d};q)_N}  \\
\!\!\!\! &=&\!\!\!\!\!\!\!
\prod_{1\leq j\leq n}
\frac{(q^{1-(n-j)g-\hat{g}_a+\hat{g}_c},
         q^{1-(n-j)g-\hat{g}_a+\hat{g}_d};q)_{-N}}
     {(q^{1-(2n-j-1)g-2\hat{g}_a},
         q^{1+(n-j)g+\hat{g}_c+\hat{g}_d};q)_{-N}}
\nonumber
\end{eqnarray}
($q\neq 1$) and
\begin{eqnarray}
\makebox[1em]{}\mathcal{N}_{1,\Lambda_N}
\!\!\!\! &=&\!\!\!\!\!\!\!
\prod_{1\leq j\leq n}\!\!
\frac{(1+(2n-j-1)g+2g_a,1-(n-j)g-g_c-g_d)_N}
     {(1+(n-j)g+g_a-g_c,1+(n-j)g+g_a-g_d)_N}  \\
\!\!\!\! &=&\!\!\!\!\!\!\!
\prod_{1\leq j\leq n}\!\!
\frac{(1-(n-j)g-\hat{g}_a+\hat{g}_c,1-(n-j)g-\hat{g}_a+\hat{g}_d)_{-N}}
     {(1-(2n-j-1)g-2\hat{g}_a,1+(n-j)g+\hat{g}_c+\hat{g}_d)_{-N}} ,
\nonumber
\end{eqnarray}
\end{subequations}
respectively.
When $n=1$ the summation formula of Theorem~\ref{finsumRDthm} reduces to
the terminating Rogers sum (cf. \cite{rog:third,gas-rah:basic})
\begin{subequations}
\begin{eqnarray}\label{rog-uni-term}
\lefteqn{
\makebox[2em]{}\sum_{\lambda \in \{ 0,\ldots ,N\} }
q^{(1-g_a-g_b-g_c-g_d)\lambda}
\Big( \frac{1-q^{2g_a+2\lambda}}
           {1-q^{2g_a}} \Bigr)
\prod_{1\leq r\leq 4}
\frac{(q^{g_r+g_a};q)_\lambda}
     {(q^{1-g_r+g_a};q)_\lambda} } &&  \\
&&={}_6\phi_5 \Bigl(
\begin{array}{c}
q^{1+g_a},-q^{1+g_a}, q^{2g_a},q^{-N},q^{g_a+g_c},q^{g_a+g_d} \\
q^{g_a}, -q^{g_a}, q^{1+N+2g_a},q^{1+g_a-g_c},q^{1+g_a-g_d}
\end{array} ;q, q^{1+N-g_c-g_d} \Bigr)  \nonumber \\
&& =
\frac{(q^{1+2g_a},q^{1-g_c-g_d};q)_N}
     {(q^{1+g_a-g_c},q^{1+g_a-g_d};q)_N}=
\frac{(q^{1-\hat{g}_a+\hat{g}_c},q^{1-\hat{g}_a+\hat{g}_d};q)_{-N}}
     {(q^{1-2\hat{g}_a},q^{1+\hat{g}_c+\hat{g}_d};q)_{-N}} \nonumber
\end{eqnarray}
for $q\neq 1$ and to the terminating Dougall sum
(cf. \cite{dou:vandermonde,gas-rah:basic})
\begin{eqnarray}\label{doug-uni-term}
\lefteqn{
\sum_{\lambda \in \{ 0,\ldots ,N\} }\:
 \Bigl( 1+\frac{\lambda}{g_a}\Bigr)
\prod_{1\leq r\leq 4}\frac{(g_r+g_a)_\lambda}
                          {(1-g_r+g_a)_\lambda} } && \\
&& =
{}_5F_4 \left(
\begin{array}{c}
1+g_a, 2g_a,-N,g_a+g_c,g_a+g_d \\
g_a,1+N+2g_a,1+g_a-g_c,1+g_a-g_d
\end{array}; 1 \right) \nonumber \\
&& =
\frac{(1+2g_a,1-g_c-g_d)_N}
     {(1+g_a-g_c,1+g_a-g_d)_N}=
\frac{(1-\hat{g}_a+\hat{g}_c,1-\hat{g}_a+\hat{g}_d)_{-N}}
     {(1-2\hat{g}_a,1+\hat{g}_c+\hat{g}_d)_{-N}} \nonumber
\end{eqnarray}
\end{subequations}
for $q=1$,
where the parameters are assumed to satisfy the
truncation condition $g_a+g_b+N=0$.

Let us conclude with an important application of Theorem~\ref{finsumRDthm}
to the theory of (basic) hypergeometric orthogonal
polynomials in several variables.
In \cite{die-sto:multivariable} a
multivariable generalization of the ($q$)-Racah polynomials
\cite{ask-wil:set,wil:some} was studied. For $q\neq 1$ the multivariable
polynomials of interest amount to the $BC$ type
Askey-Wilson polynomials due (for special parameters)
to Macdonald \cite{mac:orthogonal}
and (for general parameters) to Koornwinder \cite{koo:askey-wilson},
whereas the degenerate case $q=1$ corresponds
to the multivariable Wilson polynomials considered in \cite{die:properties}.
In \cite{mac:orthogonal,koo:askey-wilson,die:properties}
the parameter domains were chosen such that
the relevant orthogonality measures are continuous.
The paper \cite{die-sto:multivariable} focuses on a different
parameter regime for which the orthogonality measure becomes purely
discrete and finitely supported on the grid points $\rho_a+\lambda$,
$\lambda\in\Lambda_N$, with weights that are given by
$\Delta_q(\lambda )$ \eqref{rogterm}, \eqref{douterm}.
Among other things, \cite{die-sto:multivariable} presents
product formulas for the (squared) norms of the multivariable
($q$-)Racah polynomials (with respect to the corresponding
discrete inner product) in terms of the (squared) norm of the
unit polynomial. The sum in Theorem~\ref{finsumRDthm} provides
a product formula for the latter norm (i.e. the squared norm
of the unit polynomial) and thus
completes the solution of the orthonormalization problem
for the multivariable ($q$)-Racah polynomials in product form.
For $n=1$ the resulting product formulas reduce to the norm formulas for
the one-variable ($q$-)Racah polynomials presented by
Askey and Wilson \cite{ask-wil:set,wil:some}.

\begin{remark}
Multidimensional analogues of the Dougall ${}_5F_4$ sum
\eqref{doug-uni}, \eqref{doug-uni-term}
and the Rogers ${}_6\phi_5$ sum \eqref{rog-uni}, \eqref{rog-uni-term}
different from those considered in Theorem~\ref{sumRDthm} and
Theorem~\ref{finsumRDthm}
can be found in \cite{hol:summation} (${}_5F_4$ type) and
\cite{mil:q-analog,mil:basic,lil-mil:bailey,sch:summation}
(${}_6\phi_5$ type).
\end{remark}

\section{Proof of Theorem~\ref{sumMthm}}\label{sec4}
We will prove Theorem~\ref{sumMthm} by deducing a recurrence
relation for the Macdonald type sum \eqref{sumM} using a technique
due to Gustafson, who applied it to evaluate
a Selberg type multivariable generalization of
the Askey-Wilson integral \cite{gus:generalization}.
Essential ingredient in the derivation of the
recurrence relation is a multiple summation formula
taken from \cite{gus:macdonald,gus:some}.
Let the function $\Delta_q^{\text{G}}(\mathbf{x})$
be given for $0<q<1$ by
\begin{eqnarray}\label{termGq}
\Delta_q^{\text{G}}(\mathbf{x})&=&
\prod_{1\leq j<k\leq n}
\frac{1}{(q^{1+x_j+x_k},q^{1+x_j-x_k},q^{1-x_j+x_k},q^{1-x_j-x_k};q)_\infty} \\
&& \times \prod_{1\leq j\leq n}
\frac{\prod_{r=1}^{2n+2} (q^{1+g_r+x_j},  q^{1+g_r-x_j};q)_\infty}
     { (q^{1+2x_j}, q^{1-2x_j};q)_\infty } \nonumber
\end{eqnarray}
and for $q=1$ by
\begin{eqnarray}\label{termG}
\makebox[2em]{}\Delta_1^{\text{G}}(\mathbf{x})\!\!\!\!\!&=& \!\!\!\!\!\!\!
\prod_{1\leq j<k\leq n}\!\!\!
\Gamma (1+x_j+x_k)
\Gamma (1+x_j-x_k)
\Gamma (1-x_j+x_k)
\Gamma (1-x_j-x_k) \\
&& \times \prod_{1\leq j\leq n}
\frac{\Gamma (1+2x_j) \Gamma (1-2x_j)}
     {\prod_{r=1}^{2n+2} \Gamma (1+g_r+x_j) \Gamma (1+g_r-x_j)}  .\nonumber
\end{eqnarray}
Then one has for parameters satisfying
the condition $\text{Re}(1+g_1+\cdots +g_{2n+2})>0$
that \cite{gus:macdonald,gus:some}
\begin{eqnarray}\label{sumG}
\lefteqn{
\sum_{\lambda \in \mathbb{Z}^n}
\Delta_q^{\text{G}}(\mathbf{z}+\lambda ) } && \\
&& = \left\{ \begin{array}{ll}
{\displaystyle
\frac{ (q;q)_\infty^n \prod_{1\leq r < s\leq 2n+2} (q^{1+g_r+g_s};q)_\infty}
     {(q^{1+g_1+\cdots +g_{2n+2}};q)_\infty} }
& \text{for}\;\; 0<q<1 \\ [2ex]
{\displaystyle
\frac{\Gamma (1+g_1+\cdots +g_{2n+2})}
     {\prod_{1\leq r < s\leq 2n+2} \Gamma (1+g_r+g_s)}}
 & \text{for}\;\; q=1
\end{array} \right. \nonumber
\end{eqnarray}
(where $\mathbf{z}\in \mathbb{C}^n$ is assumed to satisfy
the genericity condition that the combinations
$z_j+z_j$, $z_j-z_k$ ($1\leq j<k\leq n$) and $2z_j$ ($1\leq j\leq n$)
are nonzero modulo $\mathcal{\Omega}_q$ \eqref{lattice}).
Moreover, the series on the l.h.s.
of \eqref{sumG} converges in absolute value (this may be verified in a
similar manner as was done in the appendix at the end of the paper for the
Macdonald type series on the l.h.s. of \eqref{sumM}).
For $n=1$, Gustafson's sum in \eqref{sumG} reduces to the
summation formula \eqref{bai-sum}, \eqref{doug-sum}. In other words,
the Gustafson sum \eqref{sumG} amounts to a
multidimensional generalization of the (in essence) Bailey/Dougall
sum \eqref{bai-sum}/\eqref{doug-sum} that is
different from the Macdonald type generalization
described by Theorem~\ref{sumMthm} (cf. the remark below).

Let us abbreviate the Macdonald type sum
$\sum_{\lambda \in \mathbb{Z}^n}
 1/(\mathcal{C}_{+,q} (\mathbf{z}+\lambda )\,
    \mathcal{C}_{+,q} (-\mathbf{z}-\lambda ) )$
by $S_n(g,g_r;q)$ and consider the double sum (for $0<q<1$)
\begin{eqnarray}\label{double}
\makebox[2em]{}\sum_{\lambda \in \mathbb{Z}^n}
\sum_{\mu \in \mathbb{Z}^{n-1}}\!\!\!\!
&& \!\!\!\!\prod_{1\leq j<k\leq n}
\Bigl( \frac{1}
{(q^{1+z_j+z_k+\lambda_j+\lambda_k},
  q^{1+z_j-z_k+\lambda_j-\lambda_k};q)_\infty} \\
&&\makebox[4em]{} \times
\frac{1}
 {(q^{1-z_j+z_k-\lambda_j+\lambda_k},
  q^{1-z_j-z_k-\lambda_j-\lambda_k};q)_\infty}\Bigr)\nonumber \\
&&  \!\!\!\! \times \prod_{j=1}^n
\frac{\prod_{r=1}^{4} (q^{1+g_r+z_j+\lambda_j},
                       q^{1+g_r-z_j-\lambda_j};q)_\infty}
     { (q^{1+2z_j+2\lambda_j}, q^{1-2z_j-2\lambda_j};q)_\infty }\nonumber \\
&&  \!\!\!\! \times
\prod_{1\leq j<k\leq n-1}
\Bigl(
\frac{1}
{(q^{1+z_j+z_k+\mu_j+\mu_k},
  q^{1+z_j-z_k+\mu_j-\mu_k};q)_\infty} \nonumber \\
&&\makebox[4em]{} \times
 \frac{1}{(q^{1-z_j+z_k-\mu_j+\mu_k},
           q^{1-z_j-z_k-\mu_j-\mu_k};q)_\infty} \Bigr) \nonumber \\
&& \!\!\!\! \times \prod_{k=1}^{n-1}
\frac{1}
     {(q^{1+2z_k+2\mu_k},q^{1-2z_k-2\mu_k};q)_\infty} \nonumber \\
&& \!\!\!\! \times \prod_{j=1}^n\prod_{k=1}^{n-1}
\Bigl(
(q^{1+g/2+z_j+z_k +\lambda_j+\mu_k},
       q^{1+g/2+z_j-z_k +\lambda_j-\mu_k};q)_\infty \nonumber \\
&&\makebox[4em]{} \times
(q^{1+g/2-z_j+z_k -\lambda_j+\mu_k},
          q^{1+g/2-z_j-z_k -\lambda_j-\mu_k};q)_\infty \Bigr) . \nonumber
\end{eqnarray}
We can express this
double sum in terms of $S_n(g,g_r;q)$ by
evaluating the inner sum (in $\mu$) with
the aid of Gustafson's formula \eqref{sumG}.
Alternatively, we may also reverse the order
of summation in $\lambda$ and $\mu$ (using the convergence
in absolute value)
and instead apply Gustafson's formula \eqref{sumG}
to the (now inner) sum over $\lambda$; this yields
an expression for the double sum in terms of
$S_{n-1}(g,g_r+g/2;q)$.
Comparing the two expressions for the double
sum thus obtained entails the recurrence relation
\begin{equation}
S_n(g,g_r;q) =
\frac{(q,q^{1+ng};q)_\infty \prod_{1\leq r<s\leq 4}(q^{1+g_r+g_s};q)_\infty}
     {(q^{1+g},q^{1+(n-1)g+g_1+g_2+g_3+g_4};q)_\infty}
S_{n-1}(g,g_r+g/2;q) ,
\end{equation}
which by induction on $n$ starting from the known value
for $n=1$ (taken from the r.h.s. of \eqref{bai-sum}) produces
\begin{equation}\label{rhsGq}
S_n(g,g_r;q)=\prod_{1\leq j\leq n}
\frac{ (q, q^{1+jg};q)_\infty
       \prod_{1\leq r < s\leq 4} (q^{1+(n-j)g+g_r+g_s};q)_\infty}
     {(q^{1+g},q^{1+(2n-j-1)g+g_1+g_2+g_3+g_4} ;q)_\infty} .
\end{equation}
For $q=1$ the derivation is quite similar and leads us starting from
the $q=1$ degeneration of the double sum \eqref{double} via
Gustafson's sum \eqref{sumG} to the recurrence relation
\begin{eqnarray}
\lefteqn{S_n(g,g_r;1) =} && \\
&& \frac{ \Gamma (1 +g)\, \Gamma (1+ (n-1)g +g_1+g_2+g_3+g_4)}
     { \Gamma (1+ ng)\, \prod_{1\leq r < s\leq 4} \Gamma (1+ g_r +g_s )}\:
S_{n-1}(g,g_r+g/2;1), \nonumber
\end{eqnarray}
which---when combined with the know value for $n=1$
(from the r.h.s. of \eqref{doug-sum})---entails
\begin{equation}\label{rhsG}
S_n(g,g_r;1) =
\prod_{1\leq j\leq n}
\frac{\Gamma (1+g)\,\Gamma(1+(2n-j-1)g+g_1+g_2+g_3+g_4)}
     {\Gamma (1+jg)\, \prod_{1\leq r < s\leq 4} \Gamma (1+(n-j)g+g_r+g_s)} .
\end{equation}
The theorem now follows from the observation that $S_n (g,g_r;q)$
\eqref{rhsGq}, \eqref{rhsG} may be rewritten in the form
given by the r.h.s. of \eqref{sumM} (see Eq.~\eqref{rhsM}).

\begin{remark}
Division of Gustafson's sum \eqref{sumG}
from \cite[Section~2]{gus:some} by the middle term $\Delta^G_q (\mathbf{z})$
(corresponding to $\lambda =\mathbf{0}$),
turns it into a multidimensional generalization of
the Bailey and Dougall sums \eqref{bai-ser} and \eqref{doug-ser}
that was introduced in Theorem~5.1 (for $0<q<1$)
and Theorem~8.2 (for $q=1$) of \cite{gus:macdonald}.
These multiple Bailey/Dougall sums are
different from those described in Section~\ref{sec2.4} and are (for $0<q<1$)
also different from (but reminiscent of)
the multiple Bailey sums considered by Schlosser \cite{sch:summation}.
\end{remark}

\newpage
\setcounter{equation}{0}
\renewcommand{\theequation}{A.\arabic{equation}}
\section*{Appendix: Proof of absolute
convergence multiple (basic) hypergeometric series}\label{appA}
In this appendix it is shown that
the Macdonald type series \eqref{sumM},
the generalized Aomoto-Ito series \eqref{sumAI},
the multiple Bailey/Dougall series \eqref{sumBDa}, \eqref{sumBDb}
and the multiple Rogers/Dougall series \eqref{sumRD}
all converge in absolute value.
Similar convergence proofs can be found e.g. in
\cite{kan:q-selberg,ito:theta,mac:formal}.

Clearly it suffices
to prove the absolute convergence for either
one of \eqref{sumM}, \eqref{sumAI} or \eqref{sumBDa}/\eqref{sumBDb},
as the series in question are equivalent up to multiplication by overall
factors and the convergence of the Rogers/Dougall
series \eqref{sumRD} follows immediately from the convergence of
the Bailey/Dougall series \eqref{sumBDa}/\eqref{sumBDb} upon
specialization of $\mathbf{z}$ to the value $-\rho$ (and
keeping in mind that we have reflected
the parameters $g$, $g_r$ to $-g$, $-g_r$ thus giving rise to the
minus signs in the convergence condition \eqref{convcond}).
In view of the fact that the function
$1/(\mathcal{C}_{+,q}(\mathbf{x} )\, \mathcal{C}_{+,q}(-\mathbf{x}))$
is permutation-invariant and even in the components
$x_1,\ldots ,x_n$, we conclude from the representation in \eqref{sumM} that it
it is sufficient to show that the restricted sum over the dominant cone
$\Lambda$ \eqref{cone} converges absolutely
(the cone $\Lambda$ is a fundamental
domain for the action of the (Weyl) group that
permutes and flips the signs
of the components of $\lambda\in\mathbb{Z}^n$).
In terms of the Aomoto-Ito type representation of our series given by
\eqref{sumAI}
this means that it suffices to demonstrate that
\begin{equation}\label{posser}
\sum_{\lambda \in \Lambda}
\frac{1}
     {|\mathcal{C}_{+,q}(\mathbf{z}+\lambda )\,
      \mathcal{C}_{-,q}(\mathbf{z}+\lambda ) |}	 < \infty .
\end{equation}
For $0<q<1$, we see from the explicit formula for the terms in
\eqref{aomtermq} that the series on the l.h.s. of
\eqref{posser} converges provided
\begin{equation}
\sum_{\lambda \in \Lambda} | q^{\sum_{j=1}^n (1+2\hat{\rho}_j)\lambda_j}| <
\infty,
\end{equation}
which is the case when
$\text{Re}\: (1+2\hat{\rho}_j)=
 \text{Re}\: (1+2(n-j)g+g_1+g_2+g_3+g_4) > 0$ for $j=1,\ldots ,n$.
Here we have used that the other factors of the type
\begin{equation*}
|1-q^{2z_j+2\lambda_j}|,\;\;\;\;\;
\left| \frac{(q^{1+g_r+z_j+\lambda_j};q)_\infty}
            {(q^{-g_r+z_j+\lambda_j};q)_\infty}  \right|
\end{equation*}
and
\begin{equation*}
|1-q^{z_j\pm z_k+\lambda_j\pm\lambda_k}|,\;\;\;\;\;
\left| \frac{(q^{1+g+z_j\pm z_k+\lambda_j\pm\lambda_k};q)_\infty}
            {(q^{-g+z_j\pm z_k+\lambda_j\pm\lambda_k};q)_\infty}  \right|
\end{equation*}
entering the terms of \eqref{posser} (cf. \eqref{aomtermq}) are bounded
on the dominant cone $\Lambda$ \eqref{cone} (one has that
$\lim_{x\rightarrow +\infty} (aq^x;q)_\infty /(bq^x;q)_\infty =1$).

To analyze the convergence of \eqref{posser} for the degenerate
case $q=1$, we apply the gamma function asymptotics
$\Gamma (a+x) /\Gamma (b+x)=x^{a-b}( 1+O(1/x))$ for $x\rightarrow +\infty$
to the factors of the explicit formula for
$1/(\mathcal{C}_{+,1}(\mathbf{x} )\, \mathcal{C}_{-,1}(\mathbf{x}))$
in \eqref{aomterm1}.
This entails that for $q=1$ the series on the l.h.s. of
\eqref{posser} converges when
\begin{eqnarray}\label{estimate1}
&& \sum_{\lambda\in\Lambda}
\Bigl( \prod_{1\leq j<k\leq n}
(1+\lambda_j+\lambda_k)^{-2Re (g)}
(1+\lambda_j-\lambda_k)^{-2Re (g)}  \\
&& \makebox[8em]{}\times  \prod_{1\leq j\leq n} (1+\lambda_j)^{-3-2Re
(g_1+g_2+g_3+g_4)}
\Bigr) \; < \infty .\nonumber
\end{eqnarray}
(The ratio of
$1/|\mathcal{C}_{+,1}(\mathbf{z}+\lambda)\,
\mathcal{C}_{-,1}(\mathbf{z}+\lambda)|$
and the corresponding term from \eqref{estimate1} remains bounded
when $\lambda$ runs through the dominant cone $\Lambda$ \eqref{cone}.)
Similarly, the series on the l.h.s.
of \eqref{estimate1} converges when
\begin{equation}
\sum_{\lambda\in\Lambda}
\prod_{1\leq j\leq n} \; (1+\lambda_j)^{-3-(3-\epsilon_g) (n-j)Re (g)- 2Re
(g_1+g_2+g_3+g_4)}
< \infty
\end{equation}
where $\epsilon_g = 1$ if $Re (g) \geq 0$ and $\epsilon_g =-1$ if $Re (g)<0$,
which (upon combination with the standard fact that
$\sum_{m=1}^\infty m^{-\alpha}< \infty$ for $\text{Re}\; (\alpha) >1$) leads
us to the conclusion that also for $q=1$
the series on the l.h.s. of \eqref{posser} converges provided
$\text{Re}\: (1+2(n-j)g+g_1+g_2+g_3+g_4) > 0$ for $j=1,\ldots ,n$.

\section*{Acknowledgments}
The author would like to thank Prof. K. Aomoto
for providing copies of \cite{aom:product,ito:theta} and
pointing out Ref.~\cite{mac:formal}.
Thanks are also due to the referee for correcting a mistake
in one of the formulas.

\bibliographystyle{amsplain}

\end{document}